\documentclass[12pt]{amsart}
  \usepackage{amssymb}
 \newcommand{\bee}[1]{\begin{equation}\label{#1}}
 \newcommand{\ene}{\end{equation}}

\begin{document}

MSC2010 17D99

\medskip
\begin{center}
{\Large \bf ABOUT SOME VARIETIES OF LEIBNITZ ALGEBRAS}

\medskip
{\large Yu. Yu. Frolova, T. V. Skoraya}
\end{center}

The paper presents two new results concerning the varieties of
Leibnitz algebras. In the case of prime characteristic $p$ of the
base field constructed example not nilpotent variety of Leibnitz
algebras satisfying an Engel condition order $p$. In the case of
zero characteristic obtained new results concerning the space of
multilinear elements variety of left nilpotent of the class not more
than 3 Leibnitz algebras.

\vspace{5mm}

\section{Introduction}

Recall that a Leibnitz algebra is a vector space over a field with
bilinear product, which satisfies to the identity:
$$(xy)z\equiv (xz)y+x(yz).\eqno{(1)}$$ According to this identity
multiplication on an element of algebra becomes differentiation of
this algebra. On condition of performance identity of
anticommutativity $xy\equiv -yx$, Leibnitz identity is equivalent to
the Jacobi identity: $x(yz)+y(zx)+z(xy)\equiv 0$. Therefore if The
Leibnitz algebra satisfies to the identity $xx\equiv 0$, then it is
a Lie algebra. In particular, any Lie algebra is a Leibnitz algebra.
The converse is not true.  Обратное неверно. Note that probably
first the class of Leibnitz algebras was considered in
paper~\cite{Bloh} as a generalization of the notion of Lie algebras.

We will transform Leibnitz identity as follows: $x(yz)\equiv
(xy)z-(xz)y.$ According to the received identity any element of
Leibnitz algebra can be presented in linear combination of elements
in which brackets are arranged from left to right. Therefore we will
lower brackets in case of their left arrangement, i.e.
$(((x_{1}x_{2})x_{3})\dots x_{n})=x_{1}x_{2}x_{3}\dots x_{n}$. For
convenience we denote the operator of multiplication on the right,
for example, on an element $z$ a capital letter of $Z$, assuming
that $xz=xZ$. In particular, in our notation we obtain
$x{\underbrace{yy...y}_m}=xY^m$.

A set of algebras, that satisfy some fixed set of identities, is
called a variety of linear algebras.

Let the base field $\Phi$ has zero characteristic. In this case all
information on variety contains in multilinear elements of its
algebras.

Let $F(X,\textbf{V})$ be a relatively free algebra of variety
$\textbf{V}kj$ with a countable set of free generators
$X=\{x_{1},x_{2},\dots\}$ and let $P_{n}=P_{n}(\textbf{V})$ be a set
of all multilinear elements from
$x_{\sigma(1)},x_{\sigma(2)},\dots,x_{\sigma(n)}$ in
$F(X,\textbf{V})$. Note that in the future, for the sake of
presentation, we will denote generators of the relatively free
algebra as well other symbols. Let $\sigma$ be an element of
symmetric group $S_{n}$. We assume that as a result of the actions
of the left permutation $\sigma$ on element $x_{i_{1}}x_{i_{2}}\dots
x_{i_{n}}$ space $P_{n}$ we get an element
$x_{\sigma(i_{1})}x_{\sigma(i_{2})}\dots x_{\sigma(i_{n})}$. This
sets the action of $S_{n}$ on the space $P_{n}$, in consequence of
which $P_{n}$ becomes a module over the group ring of $\Phi S_{n}$.
Structure of $P_{n}$ as $\Phi S_{n}$-module plays an important role
and is actively studied for different varieties.

Recall that the standard polynomial of degree $n$ has the form:
$$St_{n}(x_{1},x_{2},\dots ,x_{n})=\sum_{q\in
S_{n}}(-1)^{q}x_{q(1)}x_{q(2)}\dots x_{q(n)},$$ where summation is
over the elements of symmetric group, and  $(-1)^{q}$ is equal to
$+1$ or $-1$ depending on the parity of the permutation $q$. Agree
variables in standard polynomial denoted by special symbols on top
(heck, wave and so on). For example, the standard polynomial of
degree $n$ from the variables $x_{1},x_{2},\dots,x_{n}$ will be
written as follows:
$St_{n}=\overline{x}_{1},\overline{x}_{2},\dots,\overline{x}_{n}$.
As a sign of this polynomial depends on parity permutations, its
variables will call skew-symmetric. Variables in different
skew-symmetric sets will be denoted by different symbols, for
example: $$\sum_{q\in S_{n}, p\in
S_{m}}(-1)^{q+p}x_{q(1)}x_{q(2)}\dots x_{q(n)}y_{p(1)}y_{p(2)}\dots
y_{p(m)}=\overline{x}_{1}\overline{x}_{2}\dots
\overline{x}_{n}\widetilde{y}_{1}\widetilde{y}_{2}\dots
\widetilde{y}_{m}.$$ Note that if the element contains the same
variables in different skew-symmetric sets, its sign depends on the
parity of the permutation implicitly, so the variables in this
element are called alternating. In the above notation we have the
equality:
$\overline{x}_{i(1)}\dots\overline{x}_{i(k)}\overline{x}_{i(k+1)}\dots\overline{x}_{i(n)}=-\overline{x}_{i(1)}\dots\overline{x}_{i(k+1)}$
$\overline{x}_{i(k)}\dots\overline{x}_{i(n)}$, i.e. variables of
alternating set can be interchanged by changing the sign of element
on the opposite. According to Leibnitz identity can be written:
$\overline{x}_{1}\overline{x}_{2}\overline{x}_{3}\overline{x}_{4}\equiv\overline{x}_{1}\overline{x}_{2}\overline{x}_{4}\overline{x}_{3}+\overline{x}_{1}\overline{x}_{2}(\overline{x}_{3}\overline{x}_{4})$.
directly obtain:
$\overline{x}_{1}\overline{x}_{2}\overline{x}_{3}\overline{x}_{4}\equiv
\frac{1}{2}\overline{x}_{1}\overline{x}_{2}(\overline{x}_{3}\overline{x}_{4})$,
and more generally:
$\overline{x}_{1}\overline{x}_{2}\dots\overline{x}_{2n+1}\equiv\frac{1}{2^{n}}\overline{x}_{1}(\overline{x}_{2}\overline{x}_{3})\dots(\overline{x}_{2n}\overline{x}_{2n+1})$.
In other words, starting from the second place, variables  of one
skew-symmetric set, standing close by, we can combine in parentheses
multiplying element on $\frac{1}{2}$ for each pair.

In the case of zero characteristic of the base field every identity
is equivalent to the system of multilinear identities, which is
obtained using the standard method of linearization \cite{malzev2}.
Here is an example of this process for the identity
$$x_{0}(xy)(xy)\equiv 0.$$
After linearization on variable $x$, we obtain:
$$x_{0}(x_{1}y)(x_{2}y)+x_{0}(x_{2}y)(x_{1}y)\equiv 0.$$
Full linearization is as follows:
$$x_{0}(x_{1}y_{1})(x_{2}y_{2})+x_{0}(x_{1}y_{2})(x_{2}y_{1})+x_{0}(x_{2}y_{1})(x_{1}y_{2})+x_{0}(x_{2}y_{2})(x_{1}y_{1})\equiv 0.$$
Arrange the linearization of the element $f$ to designate $linf$.

Recall that the Lie algebra is called engel if it satisfies $xY^m
\equiv0$. If the algebra $A$ satisfies to identity $x_1 x_2\dots
x_{c+1} \equiv0$, but not satisfies $x_1 x_2\dots x_c\equiv0$, it's
called nilpotent of the class not more than $c$. Save these
definitions in case of Leibnitz algebras.

\section{Example of not nilpotent variety of Leibnitz algebras with the condition of engel}

In the case of zero characteristic of base field E.I. Zelmanov in
paper~\cite{Z} has proved, that an engel Lie algebra is nilpotent.
Using this fact, in the paper~\cite{ya} second author proved that in
the case of zero characteristic of base field an engel Leibnitz
algebra is nilpotent.

In the case of nonzero characteristic $p$ of base field P.M. Cohn
has given an example of a non-nilpotent Lie algebra over a residue
field on a prime module $\mathbb{Z}_p,$ which satisfies identities
$(xy)(zt)\equiv 0$ and, $xY^{p+1} \equiv 0$, that is an example of a
non-nilpotent Engel metabelian Lie algebra.

Following the ideas of the article we construct a non-nilpotent
metabelian Leibnitz algebra. Let $\Phi$ be a field of prime
characteristic $p.$

{\bf Theorem.} A variety $\bf V$ of Leibnitz algebras over field
$\Phi,$ satisfies to the identities $xY^p \equiv 0$ and $x(yz)\equiv
0$ is non-nilpotent.

Proof. We prove that over the field $\Phi$ there is a non-nilpotent
Leibnitz algebra $M$, which satisfies to the identities $xY^p \equiv
0$ and $x(yz) \equiv 0.$

Let $W$ be a vector space over a field $\Phi,$ which has a basis $\{
e_f | \ f \in \Phi^\mathbb{N} \},$ where $\Phi^\mathbb{N}$ is the
set of all functions of natural argument with values in $\Phi.$ We
define on the space $W$ multiplication, assuming that the algebra
$W$ is an abelian Lie algebra, i.e. $e_{f_1}e_{f_2} \equiv 0.$ For
any natural number $m$ we denote by $\delta _m$ endomorphism of
vector space $W,$ which on the basic element $e_f$ take a value
$e_{\overline f},$ where

 $$ \overline{f} (i) =
\left\{
\begin{array}{rcl}
f(i), \mbox{ если } i\neq m,\\
f(i)+1, \mbox{ если } i=m.\\
\end{array}
\right.
$$

It's easy to verify that $\delta_m \delta_n = \delta_n \delta_m$ и
$\delta_m^p = \varepsilon,$ где $\varepsilon$ is identical
endomorphism of vector space $W.$ Let $x_i=\delta_i-\varepsilon, \
i=1,2,...,$ than
 $x_ix_j=(\delta_i-\varepsilon)(\delta_j-\varepsilon)=\delta_i\delta_j-
 \delta_i-\delta_j+\varepsilon=\delta_j\delta_i-
 \delta_j-\delta_i+\varepsilon=(\delta_j-\varepsilon)(\delta_i-\varepsilon)=x_jx_i.$

Under the commutation operation $L=\langle x_i | i \in \mathbb{N}
\rangle_{\Phi}$ --- the $\Phi$-hull of the set $\{ x_i| \ i \in
\mathbb{N} \}$ is an abelian Lie algebra, and $W$ is a right
$L$-module.

Required Leibnitz algebra is a direct sum of vector spaces $W$ and
$L,$ in which multiplication defined by the rule:
$$
(w_1+l_1)(w_2+l_2)=w_1l_2,
$$
where $w_1,w_2 \in W, \ l_1, l_2 \in L,$ а $w_1l_2$ is a result of
applying $l_2$ to the element $w_1,$ belonging to a vector space
$W.$ Denote the resulting algebra $M.$

The algebra $M$ satisfies to the identity $x(yz) \equiv 0.$ Indeed,
substituting in a verifiable identity elements $w_i+l_i \in M, \
i=1,2,3,$ we get
 $$(w_1+l_1)((w_2+l_2)(w_3+l_3))=(w_1+l_1)(w_2l_3)=0.$$

For any elements $t,y$ of algebra $M$ it is true the engel identity
$tY^p \equiv 0.$ Indeed, it's well known, that the binomial
coefficient
 $
\left (
\begin{array}{rcl}
p\\ m
\end{array}
\right) = \frac{p!}{m!(p-m)!} $ is divided into $p.$ Therefore for
any $i \in \mathbb{N}, \ f \in K^\mathbb{N}$ by the binomial formula
we obtain
$e_fX_i^p=e_f\underbrace{(\delta_i-\varepsilon)(\delta_i-\varepsilon)...(\delta_i-\varepsilon)}_p
 =e_f\underbrace{\delta_i\delta_i...\delta_i}_p-e_f\underbrace{\varepsilon\varepsilon...\varepsilon}_p=0.$
The last equality follows from the fact that
$e_f\underbrace{\delta_i\delta_i...\delta_i}_p=e_f.$ Note that the
result is true in the case $p=2,$ as
$e_fx_ix_i=e_f\delta_i\delta_i+e_f\delta_i
\delta_i=e_f\delta_i\delta_i-e_f\delta_i \delta_i=0.$ Namely,
$e_fX_i^p\equiv 0,$ for any natural $i$ and any function $f.$ Let
$y=\sum_s \alpha_sx_s+\sum_f\beta_fe_f$ be an arbitrary element $M,$
than $tY^p=t{(\sum_s \alpha_sX_s)}^p=\sum_s\alpha_s^ptX_s^p=0.$

Now we verify that $M$ is non-nilpotent Leibnitz algebra. Denote by
$f_{j_1,j_2,...,j_s},$ where $\{ j_1,j_2,...,j_s \}$ is a strictly
increasing set of natural numbers, the function of natural argument,
which takes the value 1 at points $j_1,j_2,...,j_s$ and value 0 in
the rest of the points, and by ${f_0}$ we denote the function, which
takes the zero value in all points. We prove by induction on the
number of factors the following formula
$$
e_{f_0}x_1x_2...x_m=\sum_{k=0}^m(-1)^k\sum_ { \{j_1,j_2,...,j_{m-k}
\} } e_{f_0}\delta_{j_1}\delta_{j_2}...\delta_{j_{m-k}},\eqno{(2)}$$
where $j_1<j_2<...<j_{m-k}, \ \{j_1,j_2,...,j_{m-k} \} \subset
\{1,2,...,m \}.$ For $m=1$ we obtained
$e_{f_0}x_1=e_{f_0}\delta_1-e_{f_0}=e_{f_1}-e_{f_0}$ and formula (2)
is true. Assume that the desired equality holds for $m-1,$ i.e.
$$e_{f_0}x_1x_2...x_{m-1}=\sum_{k=0}^{m-1} (-1)^k\sum_ { \{j_1,j_2,...,j_{m-1-k} \} } e_{f_0}\delta_{j_1}\delta_{j_2}...\delta_{j_{m-1-k}}.$$
We prove that the equality (2) is true for $m$. Multiply both sides
of last equation by $x_m,$ we obtain
$$
e_{f_0}x_1x_2...x_m=
\left( \sum_{k=0}^{m-1} (-1)^k\sum_ { \{j_1,j_2,...,j_{m-1-k} \} } e_{f_0}\delta_{j_1}\delta_{j_2}...\delta_{j_{m-1-k}} \right)x_m =
$$
$$
\sum_{k=0}^m(-1)^k\sum_ { \{j_1,j_2,...,j_{m-k} \} } e_{f_0}\delta_{j_1}\delta_{j_2}...\delta_{j_{m-k}}
$$

Thus, the formula (2) is true for any natural $m.$

From the equality
$e_{f_0}\delta_{j_1}\delta_{j_2}...\delta_{j_{m-k}}=e_{f_{j_1,j_2,...,j_{m-k}}}$
and the formula (2) we obtained that
$$
e_{f_0}x_1x_2...x_m =\sum_{k=0}^{m} \sum_{ \{ j_1,j_2,...,j_{m-k} \} } (-1)^k e_{f_{j_1,j_2,...,j_{m-k}}}.
$$
The last element is a linear combination of various basic elements
with coefficients 1 or -1, therefore, it is different from zero for
any natural$m.$ The theorem is proved.

{\bf Note.} If $M$ is a Lie algebra, than the identity $x(yz) \equiv
0$ is equal to the identity $yzx \equiv 0,$ and the algebra $M$ is
nilpotent.

\section{The structure of multilinear part of variety of left nilpotent of the class not more
than 3 Leibnitz algebras}

In this section we will consider the variety of Leibnitz algebras
satisfying the identity
$$x(y(zt))\equiv 0,\eqno{(3)}$$ which we will denote by
$_{3}\textbf{N}$.

The variety $_{3}\textbf{N}$ has been studied by some authors (see,
e.g., \cite{abanina}, \cite{Mishchenko}, \cite{abanina-mishchenko}).
In particular, it was proved that the variety $_{3}\textbf{N}$ has
only subvariety $\widetilde{V}_{2}$ and $\widetilde{V}_{3}$ with
almost polynomial growth (see \cite{Mishchenko}); was to find a
basis of multilinear part of variety $_{3}\textbf{N}$ (see
\cite{abanina}); were found multiplicity $m_{\lambda}$ in the
decomposition $\chi(_{3}\textbf{N})$ (see
\cite{abanina-mishchenko}).

Here is an example of the Leibnitz algebra lying in the manifold
$_{3}\textbf{N}$. Let $T_{s}=K[t_{1},\dots,t_{s}]$ be a ring of
polynomial from variable $t_{1},\dots,t_{s}$. Consider the
Heisenberg algebra $H_{s}$ with basis
$\{a_{1},\dots,a_{s},b_{1},\dots,b_{s},c\}$ and multiplying
$a_{i}b_{j}=\delta_{ij}c$, where $\delta_{ij}$ is Kronecker delta,
the product of other basic elements is equal to zero. It's well
known (see, e.g. \cite{abanina}), that the algebra $H_{s}$ is
nilpotent of the class 2 Lie algebra. Transform the polynomial ring
$T_{s}$ in the right module of algebra $H_{s}$, in which the basic
elements of $H_{s}$ act right on the polynomial $f$ from $T_{s}$ as
follows:
$$fa_{i}=f'_{i}, fb_{i}=t_{i}f, fc=f,$$ where $f'_{i}$ is
a partial derivative of polynomial $f$ in the variable $t_{i}$. It's
easy to prove that the direct sum of vector space $H_{s}$ and
$T_{s}$ with multiplication:
$$(x+f)(y+g)=xy+fy,$$ where $x,y$ are from $H_{s}$; $f,g$ are from $T_{s}$,
is a Leibnitz algebra. Denote this algebra by symbol $H^{s}$. The
resulting algebra $H^{s}$ belongs to the variety $_{3}\textbf{N}$
for any $s$. Prove that the identity (3) holds in $H^{s}$:
$$(x_{1}+f_{1})((x_{2}+f_{2})((x_{3}+f_{3})(x_{4}+f_{4})))=x_{1}(x_{2}(x_{3}x_{4})+f_{1}(x_{2}(x_{3}x_{4}))=0.$$
This equality is true by the nilpotency of class 2 of algebra
$H_{s}$.

We define the general form of the elements that are not equal to
zero. Since the Heisenberg algebra $H_{s}$ is nilpotent of class 2,
then the product of any tree elements is equal to zero. Therefore in
algebra $H^{s}$ all of elements of degree 3 and more must contain at
least one polynomial from $T_{s}$. As $T_{s}$ is a right module of
algebra $H_{s}$, then the elements of form $x_{1}f$ are zero.
Besides if the element of algebra $H^{s}$ has two polynomials from
$T_{s}$, then  he also is zero according to the definition of the
algebra $H^{s}$. Therefore all nonzero elements of algebra $H^{s}$
of degree more than two must to have exactly one polynomial at the
first place. If the element of algebra $H^{s}$ has the polynomial
outside of the first skew-symmetric set, then it can be written as
the sum of terms, each of which not contains a polynomial in the
first place and so it's equal to zero. Therefore nonzero elements of
algebra $H^{s}$ have the polynomial $f$ in the first skew-symmetric
set.

As previously mentioned, the space of multilinear elements of degree
$n$ of some variety of Leibnitz algebras is a direct sum of
irreducible submodules, corresponding to all Young diagrams, which
contain $n$ cells; moreover, two module are isomorphic if and only
if they meet the same diagram. In the paper \cite{abanina} it was
proved, that the number of isomorphic terms in specified sum for the
space $P_{n}(_{3}\textbf{N})$ is equal to the number of corner cell
corresponding Young diagram. Therefore we will consider the diagrams
with fixed number of corner cell.

Given a Young diagram with $n$ cells, which has $k$ corner cells,
i.e. it responds to the partition $\lambda
=(n_{1}^{m_{1}},n_{2}^{m_{2}},\dots,n_{k}^{m_{k}})$, where
$n_{1}>n_{2}>\dots >n_{k}>0$ and
$n_{1}m_{1}+n_{2}m_{2}+\dots+n_{k}m_{k}=n$, that is the diagram of
form: \vskip 0.1in
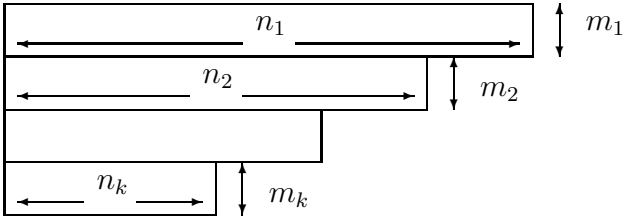
\begin{picture}(400,40)
\put(0,40){\line(1,0){200}} \put(0,20){\line(1,0){200}}
\put(0,0){\line(1,0){160}} \put(0,-20){\line(1,0){120}}
\put(0,-40){\line(1,0){80}} \put(200,40){\line(0,-1){20}}
\put(160,20){\line(0,-1){20}} \put(120,0){\line(0,-1){20}}
\put(80,-20){\line(0,-1){20}} \put(0,40){\line(0,-1){80}}
\put(95,30){$n_{1}$} \put(110,25){\vector(1,0){85}}
\put(90,25){\vector(-1,0){85}} \put(75,10){$n_{2}$}
\put(90,5){\vector(1,0){65}} \put(70,5){\vector(-1,0){65}}
\put(35,-30){$n_{k}$} \put(50,-35){\vector(1,0){25}}
\put(30,-35){\vector(-1,0){25}} \put(220,30){$m_{1}$}
\put(210,30){\vector(0,1){10}} \put(210,30){\vector(0,-1){10}}
\put(180,5){$m_{2}$} \put(170,10){\vector(0,1){10}}
\put(170,10){\vector(0,-1){10}} \put(100,-35){$m_{k}$}
\put(90,-30){\vector(0,1){10}} \put(90,-30){\vector(0,-1){10}}
\end{picture}

\vskip 0.6in

Consider the special case of diagrams of this form. Let $n=21$ and
$\lambda=(6,6,4,4,1)$. Then responding diagram has form:

\vskip 0.1in
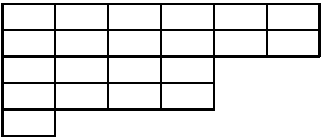
\begin{picture}(400,40)
\put(0,40){\line(1,0){120}} \put(0,30){\line(1,0){120}}
\put(0,20){\line(1,0){120}} \put(0,10){\line(1,0){80}}
\put(0,0){\line(1,0){80}} \put(0,-10){\line(1,0){20}}
\put(120,40){\line(0,-1){20}} \put(100,40){\line(0,-1){20}}
\put(80,40){\line(0,-1){40}} \put(60,40){\line(0,-1){40}}
\put(40,40){\line(0,-1){40}} \put(20,40){\line(0,-1){50}}
\put(0,40){\line(0,-1){50}}
\end{picture}

\vskip 0.1in

To this diagram respond next elements:
$$g_{1}=\overline{x}_{1}\overline{x}_{2}\overline{x}_{3}\overline{x}_{4}\overline{x}_{5}\widetilde{x}_{1}\widetilde{x}_{2}\widetilde{x}_{3}\widetilde{x}_{4}\widehat{x}_{1}\widehat{x}_{2}\widehat{x}_{3}\widehat{x}_{4}\overline{\overline{x}}_{1}\overline{\overline{x}}_{2}\overline{\overline{x}}_{3}\overline{\overline{x}}_{4}\widetilde{\widetilde{x}}_{1}\widetilde{\widetilde{x}}_{2}\widehat{\widehat{x}}_{1}\widehat{\widehat{x}}_{2},$$
$$g_{2}=\overline{x}_{1}\overline{x}_{2}\overline{x}_{3}\overline{x}_{4}\widetilde{x}_{1}\widetilde{x}_{2}\widetilde{x}_{3}\widetilde{x}_{4}\widetilde{x}_{5}\widehat{x}_{1}\widehat{x}_{2}\widehat{x}_{3}\widehat{x}_{4}\overline{\overline{x}}_{1}\overline{\overline{x}}_{2}\overline{\overline{x}}_{3}\overline{\overline{x}}_{4}\widetilde{\widetilde{x}}_{1}\widetilde{\widetilde{x}}_{2}\widehat{\widehat{x}}_{1}\widehat{\widehat{x}}_{2},$$
$$g_{3}=\overline{x}_{1}\overline{x}_{2}\widetilde{x}_{1}\widetilde{x}_{2}\widetilde{x}_{3}\widetilde{x}_{4}\widetilde{x}_{5}\widehat{x}_{1}\widehat{x}_{2}\widehat{x}_{3}\widehat{x}_{4}\overline{\overline{x}}_{1}\overline{\overline{x}}_{2}\overline{\overline{x}}_{3}\overline{\overline{x}}_{4}\widetilde{\widetilde{x}}_{1}\widetilde{\widetilde{x}}_{2}\widetilde{\widetilde{x}}_{3}\widetilde{\widetilde{x}}_{4}\widehat{\widehat{x}}_{1}\widehat{\widehat{x}}_{2}.$$

Using the notation introduced earlier for the operators, agree
standard polynomial of operators $X_{1},X_{2},\dots,X_{m}$ denoted
by
$\overline{St}_{m}=\overline{St}_{m}(X_{1},\dots,X_{m})=\overline{X}_{p(1)}\dots
\overline{X}_{p(n)}$. Note that the standard polynomials, containing
different alternating set of variables, we will write with different
upper symbols. In this designation it also true the generalization
to the case of degree of standard polynomial from operators. Then
the elements $g_{1},g_{2}$ and $g_{3}$ we can write as follows:
$$g_{1}=\overline{x}_{1}\overline{x}_{2}\overline{x}_{3}\overline{x}_{4}\overline{x}_{5}\widetilde{St}_{4}^{3}\widehat{St}_{2}^{2},$$
$$g_{2}=\overline{x}_{1}\overline{x}_{2}\overline{x}_{3}\overline{x}_{4}\widetilde{St}_{5}\widehat{St}_{4}^{2}\overline{\overline{St}}_{2}^{2},$$
$$g_{3}=\overline{x}_{1}\overline{x}_{2}\widetilde{St}_{5}\widehat{St}_{4}^{3}\overline{\overline{St}}_{2}.$$

We will return to the general case. Let as before the Young diagram
respond to the partition $\lambda
=(n_{1}^{m_{1}},n_{2}^{m_{2}},\dots,n_{k}^{m_{k}})$, where
$n_{1}>n_{2}>\dots >n_{k}>0$ and
$n_{1}m_{1}+n_{2}m_{2}+\dots+n_{k}m_{k}=n$. Then to it corresponds
the elements of following form:
$$g_{1}=(\overline{x}_{1}\overline{x}_{2}\dots\overline{x}_{d_{k}})\widetilde{St}_{d_{k}}^{n_{k}-1}\widehat{St}_{d_{k-1}}^{n_{k-1}-n_{k}}\dots \overline{\overline{St}}_{d_{2}}^{n_{3}-n_{2}}\widetilde{\widetilde{St}}_{d_{1}}^{n_{2}-n_{1}},$$
$$g_{2}=(\overline{x}_{1}\overline{x}_{2}\dots\overline{x}_{d_{k-1}})\widetilde{St}_{d_{k}}^{n_{k}}\widehat{St}_{d_{k-1}}^{n_{k-1}-n_{k}-1}\dots \overline{\overline{St}}_{d_{2}}^{n_{3}-n_{2}}\widetilde{\widetilde{St}}_{d_{1}}^{n_{2}-n_{1}},$$
$$\dots$$
$$g_{k}=(\overline{x}_{1}\overline{x}_{2}\dots\overline{x}_{d_{1}})\widetilde{St}_{d_{k}}^{n_{k}}\widehat{St}_{d_{k-1}}^{n_{k-1}-n_{k}}\dots \overline{\overline{St}}_{d_{2}}^{n_{3}-n_{2}}\widetilde{\widetilde{St}}_{d_{1}}^{n_{2}-n_{1}-1},$$
where $d_{j}=\sum_{i=1}^{j}m_{i}$, $j=1,...,k$.

It is known (see.\cite{giambruno-zaicev}, chapter 2.4, p.54), that
the linearization of any element $g_{m}(\lambda)$, where
$m=1,\dots,k$, generates an irreducible module
$W_{m}(\lambda)=KS_{n}(lin g_{m}(\lambda))$, corresponding to the
fixed partition $\lambda$.

The main result of this section is the prove of the fact, that the
elements $g_{1},g_{2},...,g_{k}$ generate the irreducible
$S_{n}$-modules that provide the direct sum for the multilinear part
$P_{n}$.

\textbf{Theorem.}For any natural number $n\geq 1$ it's true the
equality
$$P_{n}=\bigoplus_{\lambda\vdash
n}\bigoplus_{r=1}^{k(\lambda)}W_{r}(\lambda). \eqno{(4)}
$$

Proof. In the paper \cite{mishchenko-zaicev} it was proved that the
number of isomorphic irreducible submodules in the formula (4) is
equal to the maximum number of linearly independent elements of form
$g_{i}$, $i=1,\dots,k$, which generate these submodules. So for the
proof of theorem it suffices to prove the linear independence of
elements $g_{1},\dots,g_{k}$. Assume the contrary. Suppose that
there is a linear relationship $$
\alpha_{1}g_{1}+\alpha_{2}g_{2}+\dots+\alpha_{k}g_{k}=0, \eqno{(5)}
$$
where At least one of $\alpha_{i}$, $i=1,\dots,k$ is non-zero.
Suppose that $l$ is the smallest number such that $\alpha_{l}\neq 0$
and $\alpha_{q}=0$, if $q<l$. Now let $d_{j}=2p_{j}+\varepsilon_{j}$
where $\varepsilon_{j}=1$, if $d_{j}$ is odd and
$\varepsilon_{j}=0$, if $d_{j}$ is even. Will carry out the
substitution variables in place of elements $g_{l},\dots,g_{k}$ of
elements from algebra $H^{s}$ as follows: if $\varepsilon_{j}=1$,
that $x_{d_{j}}=a_{p_{j}+1}$ ($j\neq k$), $x_{d_{k}}=f+a_{p_{k}+1}$;
if $\varepsilon_{j}=0$, that $x_{d_{j}}=b_{p_{j}}$ ($j\neq k$),
$x_{d_{k}}=f+b_{p_{k}}$. Elements resulting from the substitution
from $g_{l},\dots,g_{k}$, we will denote by $v_{l},\dots,v_{k}$
respectively. Then the elements $v_{l+1},\dots,v_{k}$ are zero by
definition of algebra $H^{s}$, as well as of their structure can be
seen, that they do not contain a polynomial $f$ in the first place.

We will consider now following private type of an element
$$g_{l}=\overline{x}_{1}\overline{x}_{2}\overline{x}_{3}\overline{x}_{4}\overline{x}_{5}\overline{x}_{6}\overline{\overline{x}}_{1}\overline{\overline{x}}_{2}\overline{\overline{x}}_{3}\overline{\overline{x}}_{4}\overline{\overline{x}}_{5}\overline{\overline{x}}_{6}\widetilde{x}_{1}\widetilde{x}_{2}\widetilde{x}_{3}\widetilde{x}_{4}\widetilde{\widetilde{x}}_{1}\widetilde{\widetilde{x}}_{2}\widetilde{\widetilde{x}}_{3}.$$
After substituting its variables, we obtain an element:
$$v_{l}=\overline{a}_{1}\overline{b}_{1}\overline{a}_{2}\overline{b}_{2}\overline{a}_{3}\overline{f+b_{3}}\overline{\overline{a}}_{1}\overline{\overline{b}}_{1}\overline{\overline{a}}_{2}\overline{\overline{b}}_{2}\overline{\overline{a}}_{3}\overline{\overline{f+b_{3}}}\widetilde{a}_{1}\widetilde{b}_{1}\widetilde{a}_{2}\widetilde{b}_{2}\widetilde{\widetilde{a}}_{1}\widetilde{\widetilde{b}}_{1}\widetilde{\widetilde{a}}_{2}\equiv$$
$$\equiv (-1)^{5}f\overline{a}_{1}\overline{b}_{1}\overline{a}_{2}\overline{b}_{2}\overline{a}_{3}\overline{\overline{a}}_{1}\overline{\overline{b}}_{1}\overline{\overline{a}}_{2}\overline{\overline{b}}_{2}\overline{\overline{a}}_{3}\overline{\overline{b}}_{3}\widetilde{a}_{1}\widetilde{b}_{1}\widetilde{a}_{2}\widetilde{b}_{2}\widetilde{\widetilde{a}}_{1}\widetilde{\widetilde{b}}_{1}\widetilde{\widetilde{a}}_{2}\equiv$$
 $$\equiv \frac{-1}{2^{8}}f(\overline{a}_{1}\overline{b}_{1})(\overline{a}_{2}\overline{b}_{2})\overline{a}_{3}(\overline{\overline{a}}_{1}\overline{\overline{b}}_{1})(\overline{\overline{a}}_{2}\overline{\overline{b}}_{2})(\overline{\overline{a}}_{3}\overline{\overline{b}}_{3})(\widetilde{a}_{1}\widetilde{b}_{1})(\widetilde{a}_{2}\widetilde{b}_{2})(\widetilde{\widetilde{a}}_{1}\widetilde{\widetilde{b}}_{1})\widetilde{\widetilde{a}}_{2}\equiv$$
 $$\equiv \frac{-1}{2^{8}}fcccccccca_{3}a_{2}\equiv \frac{-1}{2^{8}}fa_{3}a_{2}.$$

We return to the general case. According to the arguments given
earlier, in the element $v_{l}$ are non-zero only those terms, that
contain exactly one polynomial $f$ in the first place; that is among
the first $n_{k}$ of a skew-symmetric sets (which include one of the
sums $f+a_{l_{k}+1}$ or $f+b_{p_{k}}$), only the first set contains
a polynomial $f$, and the rest contains the second term
$a_{p_{k}+1}$ or $b_{p_{k}}$. In the first skew-symmetric set the
polynomial $f$ is initially last. So that it was on the first place,
it's necessary $d_{k}-1$ times rearrange it with elements to the
left of it. So the element $v_{l}$  will receive a coefficient
$(-1)^{d_{k}-1}$. In addition the element $v_{l}$ after polynomial
$f$ or elements $a_{i}b_{i}$, standing near, or elements $a_{i}$
that stand alone. The pair of elements $a_{i}b_{i}$ combine in
brackets, so that each pair gives the coefficient $\frac{1}{2}$. The
product of elements $a_{i}b_{i}$ inside the brackets is equal to
$c$, since the Kronecker delta $\delta_{ii}$ is equal to 1. Since
the multiplication polynomial $f$ by $c$ gives again a polynomial
$f$, then in element $v_{l}$ to the right of $f$ will remain only
elements $a_{i}$ that stand outside the brackets. Therefore, the
element $v_{l}$ will have the form :
$$v_{l}=\frac{(-1)^{d_{k}-1}}{2^{r}}fa_{p_{k}+\varepsilon_{k}}^{
\varepsilon_{k}(n_{k}-2)+1}a_{p_{k-1}+1}^{\varepsilon_{k-1}(n_{k-1}-n_{k})}\dots
a_{p_{1}+1}^{\varepsilon_{1}(n_{1}-n_{2})},$$ where
$r=p_{1}(n_{1}-n_{2})+p_{2}(n_{2}-n_{3})+\dots+p_{k-1}(n_{k-1}-n_{k})+p_{k}n_{k}+\varepsilon_{k}-1$.
The polynomial $f$ we can select so that the result of its
differentiation on variables $t_{1},\dots,t_{s}$ was nonzero. Thus,
the element $v_{l}$ is nonzero. Note that the elements
$v_{l+1},...,v_{k}$ in such a substitution are equal to zero,
because they does not contain a polynomial of the rings $T_{s}$ in
the first skew-symmetric set.

Let us return to the linear combination (5). By assumption
$\alpha_{1}=\alpha_{2}=\dots=\alpha_{l-1}=0$. After making this
substitution this linear combination will have the form:
$\alpha_{l}\cdot v_{l}+\alpha_{l+1}\cdot 0 +\dots+\alpha_{k}\cdot
0=0$, where $v_{l}\neq 0$. Therefore, $\alpha_{l}=0$. We received a
contradiction with the assumption. Therefore, all the coefficients
$\alpha_{1},\alpha_{2},\dots,\alpha_{k}$ in the linear combinations
(5) are equal to zero, that is, elements $g_{1},g_{2},\dots,g_{k}$
are linearly independent. The theorem is proved.

\vskip 0.2in

The authors would like to thank their supervisor S.P. Mishchenko for
the formulation of a problem, useful tips and attention to the
paper.

\end{document}